\documentclass[12pt]{article}
\usepackage[utf8]{inputenc}
\usepackage[T1]{fontenc}
\usepackage{graphicx}
\usepackage{parskip}
\usepackage[english]{babel}
\usepackage{hyperref}
\usepackage{footmisc}
\usepackage{amsmath}
\usepackage{amssymb}
\usepackage{amsthm}

\usepackage{cleveref}
\usepackage{array}
\usepackage{pdfpages}

\newcommand{\uproman}[1]{\uppercase\expandafter{\romannumeral#1}}

\usepackage{multirow}
\usepackage{hhline}
\usepackage{makecell}

\usepackage{physics}
\usepackage{mathtools}
\usepackage{wasysym}
\usepackage{amsxtra}
\usepackage{comment}
\usepackage{float}
\usepackage{lmodern}
\usepackage{verbatim}
\usepackage{amsfonts}
\usepackage{booktabs, boldline}
\usepackage{siunitx}
\usepackage{tabularx}
\newcolumntype{C}{>{\centering\arraybackslash}X}
\newcolumntype{R}{>{\raggedleft\arraybackslash}X}
\newcolumntype{L}{>{\raggedright\arraybackslash}X}
\usepackage{tabularx, caption, booktabs, makecell}
\newtheorem{conjecture}{Conjecture}

\usepackage{geometry}
\usepackage{ragged2e}
\usepackage{caption}
\usepackage{setspace}
\usepackage[
backend=biber,
style=numeric,
sorting=none,
maxnames=99
]{biblatex}
\usepackage{sectsty}

\providecommand{\keywords}[1]
{
  \small	
  \textbf{\textit{Keywords---}} #1
}

\sectionfont{\fontsize{12}{15}\selectfont}

\title{
Modular Periodicity of Random Initialized Recurrences%:\\ A Study of 
}
\author{Marc T. Pudelko}
\date{Göttingen and Föhr, 2025}

\begin{document}
\maketitle

\begin{abstract}
Classical studies of the Fibonacci sequence focus on its periodicity modulo $m$ (the Pisano periods) with canonical initialization. We investigate instead the complete periodic structure arising from all $m^2$ possible initializations in $(\mathbb{Z}/m\mathbb{Z})^2$.
We discover perfect mirror symmetry between the Fibonacci recurrence $a_n = a_{n-1} + a_{n-2}$ and its parity transform $a_n = - a_{n-1} + a_{n-2}$ and observe fractal self-similarity in the extension from prime to prime power moduli. Additionally, we classify prime moduli based on their quadratic reciprocity and demonstrate that periodic sequences exhibit weight preservation under modular extension. Furthermore, we define a minima distribution $P(n)$ governed by Lucas ratios, which satisfies the symmetric relation $P(n)=P(1-n)$.
For cyclotomic recurrences, we propose explicit counting functions for the number of distinct periods with connections to necklace enumeration. These findings imply potential connections to Viswanath's random recurrence, modular forms and L-functions.
\end{abstract}


\section{Introduction}
Linear homogeneous recurrence relations play a fundamental role in number theory and combinatorics. A recurrence of order $k$ is defined as
\begin{equation} \label{eqn:recurrence}
a_n=r_1a_{n-1}+r_2a_{n-2}+...+r_{k-1}a_{n-k+1}+r_ka_{n-k}=\sum\limits_{i=1}^kr_ia_{n-i}
\end{equation}

where the sequence is fully determined by $k$ initial values ($a_0,a_1,...,a_{k-1}$) and its characteristic polynomial 

\begin{equation}\label{eqn: polynomial}
    f(x)=x^k-r_1x^{k-1}-...-r_{k-1}x-r_k
\end{equation}

as shown in \cite{article}. 
The most celebrated example is the Fibonacci sequence $\{0, 1, 1, 2, 3, 5, 8,13, \ldots\}$, generated by $x^2-x-1$ with initial values $a_0=0, a_1=1$. This sequence grows exponentially, but when reduced modulo $m$, it becomes periodic with period length $\pi(m)$, known as the Pisano period \cite{Willrich2019PisanoPA}. Determining $\pi(m)$ explicitly remains an open problem in number theory.\\
However, the classical Pisano problem addresses only a single sequence with fixed initialization. In $(\mathbb{Z}/m\mathbb{Z})^2$, there are $m^2$ possible initial conditions --- and in general there are $m^k$ initial conditions for a recurrence of order $k$ --- clustering into a variety of irreducible and cyclically equivalent sequences. For example, the period [011011] reduces to [011], which is cyclic equivalent to [110] and [101], capturing the three initial conditions $(0,1)$,$(1,1)$ and $(1,0)$. This naturally leads to a broader question: what is the complete structure of all periodic sequences arising from a given recurrence modulo $m$? How many distinct periods exist and what are their lengths? Understanding this "period landscape" reveals algebraic structure invisible when studying only the canonical Fibonacci sequence.\\
To approach this systematically, we investigate two families of recurrences that exhibit complementary structure. The first family consists of cyclotomic polynomials $\Phi_n(x)$, defined as
\begin{equation}
    \Phi_n(x) = {\displaystyle \prod_{\substack{1\leq k \leq n \\\gcd(k,n)=1}}\left(x-e^{2i\pi \frac{k}{n}}\right)}
\end{equation}
satisfying $\displaystyle \prod_{d\mid n} \Phi_d(x) = x^n-1$. These 
polynomials are fundamental in algebra and number theory, yet their role 
in modular periodicity remains largely unexplored. Since all their roots are simple and lie on the unit circle, cyclotomic recurrences remain bounded, making them a natural starting point for classification.\\
The second family consists of exponentially growing sequences, exemplified by the Fibonacci recurrence $a_n = a_{n-1} + a_{n-2}$ (from $x^2 - x - 1$) and its parity transform $a_n = -a_{n-1} + a_{n-2}$ (from $x^2 + x - 1$). Interestingly, these two form the deterministic branches of Viswanath's random recurrence $a_n = \pm a_{n-1} + a_{n-2}$, which also exhibits exponential growth with rate approximately $1.13198824$ \cite{article2}. Moreover, they arise as sign-reversed versions of the 
cyclotomic recurrences from $\Phi_3(x) = x^2 + x + 1$ and 
$\Phi_6(x) = x^2 - x + 1$. This connection between unit-root-bounded and exponentially growing recurrences through coefficient sign reversal motivates our parallel investigation of both families.\\
This paper investigates the period landscapes of both cyclotomic and non-cyclotomic recurrences across three main sections. In Section 2, we conjecture explicit counting functions for the number of distinct periods modulo $m$ for various cyclotomic families, including $\Phi_p$, $\Phi_{2p}$, $\Phi_{p^j}$ and $x^n-1$, all supported by extensive computational verification. In Section 3, the non-cyclotomic Fibonacci and parity recurrences are explored, revealing unexpected mirror symmetry. Despite having algebraically distinct characteristic roots, these recurrences produce identical period structures for all tested moduli. We classify their prime moduli via the Legendre symbol and a parameter $\alpha$ governing both period counts and root orders in field extensions, observe fractal self-similarity at prime power moduli and conjecture weight preservation under modular extension. Furthermore, we establish that both recurrences share an identical minimum distribution governed by Lucas number ratios. In Section 4, we discuss broader implications and open questions and suggest that deterministic periodicity in finite rings has potential relations to stochastic growth in random recurrences and to classical number-theoretical objects such as modular forms and L-functions.

\section{Periodicity of Cyclotomic Recurrences}\label{section2}

This section addresses the periodicity of recurrences with characteristic polynomial being cyclotomic. For example, following equations \ref{eqn:recurrence} and  \ref{eqn: polynomial}, the corresponding recurrence for the first cyclotomic polynomial is $a_n=a_{n-1}$, for the third it is $a_n=-a_{n-1}-a_{n-2}$ and for the ninth it is $a_n=-a_{n-3}-a_{n-6}$. Notably the first cyclotomic polynomial having other coefficients than 1,0, or -1 is $\Phi_{105=3\cdot5\cdot7}$.

\begin{conjecture}[Period Count for $\Phi_p$]
\label{conj:cyclotomic_periods}
Let $p$ be a prime number and consider the linear recurrence relation of order $p-1$ defined by the characteristic polynomial $\Phi_p(x) = \sum_{k=0}^{p-1} x^k$, where $\Phi_p(x)$ is the $p$-th cyclotomic polynomial. For a given positive integer $m$, let $\#(\Phi_p, m)$ denote the number of distinct periods when the recurrence is computed modulo $m$ over all possible initial conditions $(a_0, a_1, \ldots, a_{p-2}) \in (\mathbb{Z}/m\mathbb{Z})^{p-1}$. Then the following formula holds:
\begin{equation}
\#(\Phi_p, m) = \begin{cases}
\frac{m^{p-1} - p}{p} + p & \text{if } p \mid m \\
\frac{m^{p-1} - 1}{p} + 1 & \text{if } p \nmid m
\end{cases}
\end{equation}
\textbf{Remark}: The conjecture implies that when $p$ divides $m$, there are exactly $p$ fixed points (periods of length 1) and $\frac{m^{p-1} - p}{p}$ periods of length $p$. When $p$ does not divide $m$, there is exactly one fixed point (the zero vector) and $\frac{m^{p-1} - 1}{p}$ periods of length $p$.\\
\textbf{Example}:
For $p = 5$ and $m = 10$ (where $5 \mid 10$), we have:
\[
\#(\Phi_5,10) = \frac{10^4 - 5}{5} + 5 = \frac{9995}{5} + 5 = 1999 + 5 = 2004
\]
This corresponds to 5 fixed points (namely [0],[2],[4],[6] and [8]) and 1999 periods of length 5.
\end{conjecture}

\begin{conjecture}[Period Count for $\Phi_{2p}$]
\label{conj:cyclotomic_periods2}
Let $p$ be an odd prime number and consider the linear recurrence relation of order $p-1$ defined by the characteristic polynomial $\Phi_{2p}= \sum_{k=0}^{p-1} (-x)^k$. For a given positive integer $m$, let $\#(\Phi_{2p}, m)$ denote the number of distinct periods when the recurrence is computed modulo $m$ over all possible initial conditions $(a_0, a_1, \ldots, a_{n-1}) \in (\mathbb{Z}/m\mathbb{Z})^{p-1}$. Then the following formula holds:
\begin{equation}
    \#(\Phi_{2p}, m)=\begin{cases}
    \frac{m^{p-1} - p}{2p} + 1 + \frac{p-1}{2} & \text{if } p \mid m\\
\frac{m^{p-1}- 2^{p-1}}{2p} + 1 + \frac{2^{p-1}-1}{p} & \text{if } 2 \mid m \\
\frac{m^{p-1} -2^{p-1}-p+1}{2p} + 1+ \frac{p-1}{2}+\frac{2^{p-1}-1}{p} & \text{if } p \:\text{and } 2\mid m\\
\end{cases}
\end{equation}
\textbf{Remark}:
The conjecture implies, among other things, that when $p$ and 2 divide $m$, there is one fixed point (the zero vector), $\frac{p-1}{2}$ periods of length 2, $\frac{2^{p-1}-1}{p}$ periods of length p and  $\frac{m^{p-1}-2^{p-1}-p+1}{2p}$ periods of length $2p$.\\
\textbf{Example}:
For $2p = 10$ and $m =10$ (where  $p$ and 2 devide m), we have:
\[
\#(\Phi_{10},10) = \frac{10^{4} -  2^{4}-5+1}{10} + 1 +\frac{4}{2} + \frac{2^{4}-1}{5}= \frac{9980}{10} + 6 = 998 +6 = 1004
\]
This corresponds to 1 fixed point, 2 periods of length 2 (namely [28] and [46]), 3 periods of length 5 (namely [00055], [00505] and [05555]) and 998 periods of length 10.
\end{conjecture}

\begin{conjecture}[Period Count for $\Phi_{p^j}$]
\label{conj:cyclotomic_periods3}
Let $p$ be a prime number and consider the linear recurrence relation of order $p^j-p^{j-1}$ defined by the characteristic polynomial $\Phi_{p^j}(x)=\sum_{k=0}^{p-1}x^{kp^{j-1}}$, where $\Phi_{p^j}(x)$ is the $p^j$-th cyclotomic polynomial. For a given positive integer $m$, let $\#(\Phi_{p^j}, m)$ denote the number of distinct periods when the recurrence is computed modulo $m$ over all possible initial conditions $(a_0, a_1, \ldots, a_{p^j-p^{j-1}-1}) \in (\mathbb{Z}/m\mathbb{Z})^{p^j-p^{j-1}}$. Then the following formula holds:

\begin{equation}\label{eqn: Phi}
    \#(\Phi_{p^j}, m)= \begin{cases}
        \frac{m^{p^j-p^{j-1}}- \sum^{j-1}\limits_{i=0}p^iM(p,p^i)}{p^j} + \sum^{j-1}\limits_{i=0}M(p,p^i) & \text{if } p\mid m\\
       \frac{m^{p^j-p^{j-1}}-1}{p^j}+1 & \text{if } p \nmid m
    \end{cases}
\end{equation}

, where the periodic lengths are all $p^i$ with $i \in \{0,1,2,...,j\}$ if $p\mid m$ and $i \in \{0,j\}$ if $p\nmid m$, and where

\begin{equation}\label{eqn: M}
M(m,r)=\frac{1}{r}\sum\limits_{d\mid r}\mu (d)m^{r/d}   
\end{equation}

is the number of different aperiodic $m$-ary necklaces of length $r$, as presented in \cite{meštrović2018differentclassesbinarynecklaces}.It is also the number of monic irreducible polynomials of degree r over a finite field $\mathbb{F}_q$ (see \cite{çakıroğlu2020numberirreduciblepolynomialsfinite}), with $\mu$ being the classic Möbius function. $M(m,r)$ also refers to Moureau's necklace-counting function or MacMahon's formula.\\
\textbf{Example}:
For $p^j = 3^2=9$ and $m = 12$ (where $3 \mid 12$), we have:
\[
\#(\Phi_9,12) = \frac{12^6 - (1\cdot3 +3\cdot8)}{9} +( 3 +8) = \frac{2985957}{9} + 11 = 331773 + 27 = 331784
\]
This corresponds to 3 fixed points (namely [0],[4] and [8]), 8 periods of length 3 (namely [004], [008], [044],[088], [048], [084], [448] and [884]) and 331773 periods of length 9.
\end{conjecture}

\begin{conjecture}[Period Count for $x^n-1$]
\label{conj:cyclotomic_periods4}
Consider the linear recurrence relation of order $n$ defined by the characteristic polynomial $x^n-1$. For a given positive integer $m$, let $\#(n, m)$ denote the number of distinct periods when the recurrence is computed modulo $m$ over all possible initial conditions $(a_0, a_1, \ldots, a_{n-1}) \in (\mathbb{Z}/m\mathbb{Z})^{n}$. Then the following formula holds:

\begin{equation}
    \#(n, m)=
         \frac{1}{n}\left(m^{n}- \displaystyle\sum_{\substack{r \mid n \\ r \neq n}} rM(m,r)\right) + \sum_{\substack{r \mid n \\ r \neq n}}M(m,r) 
\end{equation}
, where the periodic lengths are all $r$ that divide n, and $M(m,r)$ is defined by \autoref{eqn: M}.\\
\textbf{Example}:
For $n = 6$ and $m =4$, we have:
\[
\#(6,4) = \frac{4^6 - (1\cdot4 +2\cdot6+3\cdot20)}{6} + (4 +6+20) = \frac{4020}{6} + 30 = 670 + 30 = 700
\]
This corresponds to 4 fixed points (namely [0],[1],[2] and [3]), 6 periods of length 2 (namely [01], [02], [03],[12], [13] and [23]), 20 periods of length 3 (namely [001], [002], [003], [011], [022], [033], [012], [013], [021], [023], [031], [032], [112], [113], [221], [223], [331], [332], [123] and [132]) and 670 periods of length 6.
\end{conjecture}

\section{Periodicity of Fibonacci Recurrences}\label{section3}
We now turn our focus to the Fibonacci recurrence $a_n=a_{n-1}+ a_{n-2}$ and its parity transform $a_n=-a_{n-1}+ a_{n-2}$ with random integer initialization $(a_0, a_1) \in (\mathbb{Z}/m\mathbb{Z})^2$. The implied mirror
symmetry can be seen in  \autoref{table1} for all moduli m, where the number of periods and lengths are
similar for both recurrences, and where every period corresponds to a mirror period from the other
recurrence.

\begin{table}[H]
\centering
\renewcommand{\arraystretch}{1.3} % Increase row spacing
\captionof{table}{Periods for different moduli $m$, with a mirror drawn between the two recurrences, and where the Pisano periods are in bold}
\label{table1}
\begin{tabularx}{\textwidth}{>{\centering\arraybackslash}p{1cm} | >{\raggedleft\arraybackslash}X || >{\raggedright\arraybackslash}X}
\hlineB{3}
\makecell[c]{m} & \makecell[c]{$a_n=a_{n-1}+ a_{n-2}$} & \makecell[c]{$a_n=-a_{n-1}+ a_{n-2}$}\\
\hlineB{3}
1 & \textbf{0} & 0\\
2 & \textbf{110} ,0 & 0, 011\\
3 & \textbf{11202210} ,0 & 0, 01220211 \\
4 & 332130 ,\textbf{112310} ,220 ,0 & 0, 022, 013211, 031233 \\
5 & \textbf{11230331404432022410}, 3421, 0 & 0, 1243, 01422023440413303211 \\
\makecell[c]{6} & \makecell[r]{,22404420 ,330 ,0\\ \textbf{112352134150554314532510}} & \makecell[l]{0, 033, 02440422,\\ 015235413455051431253211} \\
\hlineB{3}
\end{tabularx}
\end{table}

One open problem in mathematics is calculating the length of the Pisano periods, $\pi(m)$, explicitly. However, classical studies focus on the single canonical sequence, neglecting the complete period landscape arising from all possible initializations. We observe that periods of length 3 emerge at all even moduli, periods of length 8 at every third modulus, and periods of length 4 and 20 at every fifth modulus. These patterns suggest a fundamental principle: the periodicity of composite moduli is completely determined by the periodicity at prime moduli via the Chinese Remainder Theorem. This motivates our systematic classification of prime moduli, which we present below.

\begin{conjecture}[Period Count for Fibonacci Recurrences modulo p]
\label{conj:fibonacci_prime_periods}
Let $p$ be a prime and let $\#(p)$ denote the number of distinct periods 
when the Fibonacci recurrence or its parity transform is computed modulo 
$p$ over all possible initial conditions $(a_0, a_1) \in (\mathbb{Z}/p\mathbb{Z})^2$. 
The count is determined by the Legendre symbol $\left(\frac{5}{p}\right)$ 
and a positive integer $\alpha$ governing the Pisano period:\\
\textbf{Class A} ($p \equiv 2, 3 \pmod{5}$): These primes satisfy 
$\left(\frac{5}{p}\right) = -1$ and have Pisano period $\pi_A(p) = \frac{2(p+1)}{\alpha}$ 
for odd $\alpha$. The number of distinct periods is:
\begin{equation}
\#_{A}(p) = \frac{\alpha}{2}(p-1)+ 1
\end{equation}
These primes exhibit two period lengths: the zero vector and non-trivial 
periods of length $\pi_A(p)$.\\
\textbf{Class B} ($p \equiv 1, 4 \pmod{5}$): These primes satisfy 
$\left(\frac{5}{p}\right) = 1$ and have Pisano period $\pi_B(p) = \frac{p-1}{\alpha}$ for any positive integer $\alpha$. They divide into two disjoint subclasses:\\
\textbf{Subclass B1} (two period lengths): The number of distinct periods is:
\begin{equation}
\#_{B1}(p) = \alpha(p+1) + 1
\end{equation}
These primes exhibit two period lengths: the zero vector and non-trivial 
periods of length $\pi_B(p)$.\\
\textbf{Subclass B2} (three period lengths): The Pisano period contains exactly one zero (OEIS A053032 for $p \geq 11$). The number of distinct periods is:
\begin{equation}
\#_{B2}(p) = \alpha(p+2) + 1
\end{equation}
These primes exhibit three period lengths: the zero vector, an intermediate 
length $\frac{\pi_B(p)}{2}$ appearing $2\alpha$ times, and the length 
$\pi_B(p)$ appearing $p\alpha$ times.\\
All primes satisfying $p \equiv 11, 19 \pmod{20}$ belong to subclass B2, 
while those satisfying $p \equiv 1, 9 \pmod{20}$ can belong to either subclass. The prime $p=5$ is excluded as a degenerate case where the 
discriminant vanishes modulo $p$.\\
\textbf{Remark}:
The parameter $\alpha$ can be characterized through multiplicative orders 
in finite field extensions. For class B primes, the roots of $x^2-x-1$ lie 
in $\mathbb{F}_p$ with equal orders $(p-1)/\alpha$ (subclass B1) or distinct 
orders $(p-1)/\alpha$ and $(p-1)/(2\alpha)$ (subclass B2). For class A 
primes, roots lie in $\mathbb{F}_{p^2}$ with order $2(p+1)/\alpha$.

\end{conjecture}

\begin{conjecture}[Self-Similarity at Prime Powers]
\label{conj:prime_power_similarity}
At each transition $p^k \to p^{k+1}$ for the Fibonacci and parity 
recurrences, all period types from $p^k$ persist at $p^{k+1}$. 
Additionally, each period of length $\ell > 1$ newly introduced at $p^k$ 
(not present at $p^{k-1}$) generates periods of length $p\ell$ with 
multiplicity scaled by factor $p$. For class B2 primes with three period 
lengths, the middle-length periods maintain constant multiplicity $2\alpha$ 
across all powers.\\
\textbf{Example}: The base structure (denoted as $\text{multiplicity} \times \text{length}$) for the B2 prime $p=19$ with $\alpha=1$ is  $\{1, 2 \times 9, 19 \times 18\}$, at $p^2 = 361$ it is $\{1, 2 \times 9, 19 \times 18, 2 \times 171, 379 \times 342\}$, and at $p^3=6859$ it is $\{1, 2 \times 9, 19 \times 18, 2 \times 171, 379 \times 342, 2 \times 3249, 7219 \times 6498\}$. 
\end{conjecture}

\begin{conjecture}[Weight Preservation of Fibonacci Recurrences]
\label{conj:weightpreservation}
Let $F_m$ denote the set of all periodic sequences arising from the Fibonacci recurrence or its parity transform with random initialization $(a_0, a_1) \in (\mathbb{Z}/m\mathbb{Z})^2$, taken modulo $m$. For any divisor $d$ of a composite modulus $m$, we define the weight of a period $p$ of length $\ell_p$ in $F_d$ as $w_d(p) = \frac{\ell_p}{d^2}$. When a period $p_d \in F_d$ is extended to $F_m$, it gives rise to one or more periods $\{p_m^{(1)}, \ldots, p_m^{(k)}\} \subset F_m$ that reduce to $p_d$ modulo $d$. We conjecture that the total weight is conserved:

\begin{equation}
    w_d(p_d) = \sum_{i=1}^{k} w_m(p_m^{(i)})
\end{equation}

\textbf{Example}:
The space $F_2$ contains periods $0$ (weight $\frac{1}{4}$) and $011$ (weight $\frac{3}{4}$). When we extend to $F_6$, the period $0$ extends to periods $0$ and $02240442$ with combined weight $\frac{1}{36} + \frac{8}{36} = \frac{9}{36} = \frac{1}{4}$, while the period $011$ extends to periods $033$ and $011235213415055431453251$ with combined weight $\frac{3}{36} + \frac{24}{36} = \frac{27}{36} = \frac{3}{4}$, preserving the original weights exactly. Similarly, the space $F_3$ contains periods $0$ (weight $\frac{1}{9}$) and $01120221$ (weight $\frac{8}{9}$). Upon extension to $F_6$, the period $0$ extends to periods $0$ and $033$ with combined weight $\frac{1}{36} + \frac{3}{36} = \frac{4}{36} = \frac{1}{9}$, while $01120221$ extends to periods $02240442$ and $011235213415055431453251$ with combined weight $\frac{8}{36} + \frac{24}{36} = \frac{32}{36} = \frac{8}{9}$, again preserving the weight distribution perfectly.
\end{conjecture}

\begin{conjecture}[Probability distribution of Fibonacci minima]
\label{conj:fibonacci_minimum}
Consider either the Fibonacci recurrence $a_n = a_{n-1} + a_{n-2}$ or its parity recurrence $a_n = -a_{n-1} + a_{n-2}$ with random integer initialization $\{a_0, a_1\} \in \mathbb{Z}^2$. Although these sequences diverge as $n \to \pm\infty$, they possess well-defined absolute minima. Let $F_n$ and $L_n$ denote the $n$-th Fibonacci and Lucas numbers respectively. Define $P(n)$ as the probability that a randomly initialized sequence has its absolute minimum at position $n$. Then:

\begin{equation}
P(n) = \begin{cases}
\frac{1}{4} & \text{if } n = 0  \\
\frac{1}{\pi}\left(\arctan\left(\frac{L_{n}}{L_{n+1}}\right) - \arctan\left(\frac{L_{n-2}}{L_{n-1}}\right)\right) & \text{if } n \geq 1, \, n \text{ odd}\\
\frac{1}{\pi}\left(\arctan\left(\frac{L_{n-2}}{L_{n-1}}\right) - \arctan\left(\frac{L_{n}}{L_{n+1}}\right)\right) & \text{if } n \geq 2, \, n \text{ even} \\
\end{cases}
\end{equation}

Moreover, the probability satisfies the symmetry relation $P(n) = P(1-n)$.\\
\textbf{Remark}:\label{rem:modular_connection}
The formula $P(n)$ is derived under the assumption of uniform 
random initialization. However, for bounded initializations 
$\{a_0, a_1\} \in [-N, N]^2 \cap \mathbb{Z}^2$, the ratio 
$r = a_0/a_1$ takes only finitely many rational values, 
imposing a finite resolution on the real line. We therefore 
expect
\begin{equation}
     P_N(n) = P(n) + \varepsilon(N), \quad \varepsilon(N) \to 0 
  \text{ as } N \to \infty.
\end{equation}
The precise rate is left for future work. Furthermore, the 
Lucas ratios appearing in the formula, together with their 
mediants $F_{n-1}/F_n$, generate matrices in 
$\mathrm{SL}(2,\mathbb{Z})$, suggesting connections to modular 
forms and continued fraction theory.

\end{conjecture}

\section{Perspectives}\label{section4}
This paper examined the modular periodicity of randomly initialized Fibonacci 
and cyclotomic recurrences, revealing previously unknown symmetries. We 
conjectured explicit counting formulas for cyclotomic recurrences and 
discovered perfect mirror symmetry between the Fibonacci and parity 
recurrences. We classified prime moduli based on the Legendre symbol and a parameter $\alpha$ governing both period counts and multiplicative orders of roots in field extensions. Additionally, we observed fractal self-similarity at prime power moduli, established weight preservation under modular extension, and derived a probability distribution for sequence minima governed by Lucas ratios.\\
The framework of random initialization naturally extends to higher-degree polynomials. For example, the two order-6 recurrences $a_n = \pm a_{n-3} + a_{n-6}$ (negations from $\Phi_9$ and $\Phi_{18}$) exhibit mirror symmetry analogous to the Fibonacci case for all tested moduli $m \leq 11$, though their period structure is more complex with up to five distinct lengths for prime moduli. Whether the classification principles established for quadratic recurrences extend to higher-order systems remains an open question.\\
Preliminary computation reveals systematic embedding of Fibonacci periods 
into cyclotomic recurrences: class A periods appear in $\Phi_{2(p+1)p^j/\alpha}$ 
modulo $p$, while class B1 periods appear in $\Phi_{(p-1)p^j/\alpha}$ modulo $p$. 
Class B2 primes exhibit more complex behavior, with some showing additional 
projections to composite moduli. These patterns suggest $\alpha$ encodes 
structural information about cyclotomic decompositions beyond what the 
Legendre symbol reveals. Investigating embeddings into other polynomial 
families --- such as Fekete polynomials, which relate to cyclotomic polynomials 
and L-functions \cite{chidambaram2025fekete_principal} --- could reveal universal 
principles governing recurrence periodicity in finite rings.\\
Finally, determining $\alpha$ efficiently remains an open problem. While 
$\alpha$ can be characterized through multiplicative orders of roots in field 
extensions --- requiring factorization of $p \pm 1$ --- or through cyclotomic 
embedding patterns, no efficiently computable form is known. Understanding 
whether $\alpha$ admits characterization through higher-order residue 
properties, Galois-theoretic invariants, spectral properties governing 
random recurrences, or connections to other polynomial families could 
illuminate the structure of Pisano periods and help resolve one of the 
longstanding open problems in the theory of Fibonacci sequences.

\singlespacing
\nocite{*}
\printbibliography[title=References]

\end{document}